\documentclass[journal]{IEEEtran}

\usepackage{cite}

\ifCLASSINFOpdf
   \usepackage[pdftex]{graphicx}
   \graphicspath{{./Images/}}
   \usepackage{epstopdf}
   \DeclareGraphicsExtensions{.pdf,.jpeg,.png}
\else
\fi

\usepackage{amsmath}

\usepackage{amssymb}

\usepackage{nomencl}
\makenomenclature
\usepackage{etoolbox}
\renewcommand\nomgroup[1]{%
	\item[\bfseries
	\ifstrequal{#1}{A}{Sets}{%
		\ifstrequal{#1}{D}{Variables}{%
			\ifstrequal{#1}{P}{Parameters}{}}{%
				\ifstrequal{#1}{I}{Initial Conditions}{}}}%
	]}
\usepackage{siunitx}
\usepackage{color}

\begin{document}
\title{Computationally Efficient Market Simulation Tool for Future Grid Scenario Analysis}

\author{Shariq~Riaz,~\IEEEmembership{Graduate Student Member,~IEEE,}
        Gregor~Verbi\v{c},~\IEEEmembership{Senior Member,~IEEE,}
        and~Archie~C.~Chapman,~\IEEEmembership{Member,~IEEE}% <-this % stops a space
\thanks{Shariq Riaz, Gregor Verbi\v{c} and Archie C. Chapman are with the School of Electrical and Information Engineering, The University of Sydney, Sydney, New South Wales, Australia. 	e-mails: (shariq.riaz, gregor.verbic, archie.chapman@sydney.edu.au).}% <-this % stops a space
\thanks{Shariq Riaz is also with the Department of Electrical Engineering, University of Engineering and Technology Lahore, Lahore, Pakistan.}% <-this % stops a space
\thanks{}}

\maketitle

\begin{abstract}
This paper proposes a computationally efficient electricity market simulation tool (MST) suitable for future grid scenario analysis. The market model is based on a unit commitment (UC) problem and takes into account the uptake of emerging technologies, like demand response, battery storage, concentrated solar thermal generation, and HVDC transmission lines. To allow for a subsequent stability assessment, the MST requires an explicit representation of the number of online generation units, which affects powers system inertia and reactive power support capability. These requirements render a full-fledged UC model computationally intractable, so we propose unit clustering, a rolling horizon approach, and constraint clipping to increase the computational efficiency. To showcase the capability of the proposed tool, we use a simplified model of the Australian National Electricity Market with different penetrations of renewable generation. The results 
%show that the number of online units resulting from the proposed tool is very close to the binary UC run over a week-long horizon, which is confirmed by the loadability and inertia analysis. That confirms the validity of the approach for long term future grid studies.
are verified by comparison to a more expressive and computationally-intensive binary UC, which confirm the validity of the approach for long term future grid studies. 
\end{abstract}


\nomenclature[A01]{$\mathcal{C}$}{Set of consumers $c$.}
\nomenclature[A02]{$\mathcal{G}$}{Set of generators $g$, $ \mathcal{G^\text{}} = \mathcal{G^\text{syn}} \cup \mathcal{G^\text{RES}} $.}
\nomenclature[A03]{$\mathcal{G^\text{syn}}$}{Set of synchronous generators, $\mathcal{G^\text{syn}} \subseteq \mathcal{G}$.}
\nomenclature[A04]{$\mathcal{G^\text{RES}}$}{Set of renewable generators, $\mathcal{G^\text{RES}} \subseteq \mathcal{G}$.}
\nomenclature[A05]{$\mathcal{G^\text{CST}}$}{Set of concentrated solar thermal generators, $\mathcal{G^\text{CST}} \subseteq\mathcal{G}^\text{syn}$.}
\nomenclature[A06]{$\mathcal{G^\text{r}}$}{Set of synchronous generators in region $r, \bigcup_{\mathcal{G^\text{r}}}  = \mathcal{G}^\text{}$.}
\nomenclature[A07]{$\mathcal{H}$}{Set of sub-horizons $h$.}
\nomenclature[A08]{$\mathcal{L}$}{Set of power lines $l$, $\mathcal{L} = \mathcal{L^\text{AC}} \cup \mathcal{L^\text{HVDC}}$.}
\nomenclature[A09]{$\mathcal{L^\text{AC}}$}{Set of AC power lines, $\mathcal{L^\text{AC}} \subseteq \mathcal{L}$.}
\nomenclature[A10]{$\mathcal{L^\text{HVDC}}$}{Set of HVDC power lines, $\mathcal{L^\text{HVDC}} \subseteq \mathcal{L}$.}
\nomenclature[A11]{$\mathcal{N}$}{Set of nodes $n$.}
\nomenclature[A12]{$\mathcal{N^\text{r}}$}{Set of nodes in region $r$.}
\nomenclature[A13]{$\mathcal{P}$}{Set of prosumers $p$.}
\nomenclature[A14]{$\mathcal{R}$}{Set of regions $r$.}
\nomenclature[A15]{$\mathcal{S}$}{Set of storage plants $s$.}
\nomenclature[A16]{$\mathcal{T}$}{Set of time slots $t$.}

\nomenclature[D01]{$s_{g,t}$}{Number of online units of generator $g$, $s_{g,t} \in \{0,1\}$ in BUC and $s_{g,t} \in  \mathbb{Z}_{+}$ in MST.}
\nomenclature[D02]{$u_{g,t} $}{Integer startup status variable of a unit of generator $g$, $u_{g,t} \in \{0,1\}$ in BUC and $u_{g,t} \in  \mathbb{Z}_{+}$ in MST.}
\nomenclature[D03]{$d_{g,t}$}{Integer shutdown status variable of a unit of generator $g$, $d_{g,t} \in \{0,1\}$ in BUC and $d_{g,t} \in  \mathbb{Z}_{+}$ in MST.}
\nomenclature[D04]{$\delta_{n,t}$}{Voltage angle at node $n$.}
\nomenclature[D05]{$p_{l,t}^{}$}{Power flow on line $l$.}
%\nomenclature[D06]{$p_{l,t}^{\text{AC}}$}{Power flow on line $l \in \mathcal{L^\text{AC}}$.}
%\nomenclature[D07]{$p_{l,t}^{\text{HVDC}}$}{Power flow on line $l \in \mathcal{L^\text{HVDC}}$.}
\nomenclature[D08]{$\Delta p_{l,t}$}{Power loss on line $l$.}
\nomenclature[D09]{$p_{g,t}$}{Power dispatch of generator $g$.}
\nomenclature[D10]{${p}_{p,t}^{\text{g}+/-}$}{Grid/feed-in power of prosumer $p$.}
\nomenclature[D11]{${p}_{s,t}$}{Power flow of storage plant $s$.}
\nomenclature[D12]{${p}_{p,t}^{\text{b}}$}{Battery power flow of prosumer $p$.}
\nomenclature[D13]{$e_{g,t}$}{Thermal energy stored in TES of generator $g \in \mathcal{G^\text{CST}} $.}
\nomenclature[D14]{$e_{s,t}$}{Energy stored in storage plant $s$.}
\nomenclature[D15]{$e_{p,t}^\text{b}$}{Battery charge state of prosumer $p$.}
	
\nomenclature[P]{${c}_g^\text{fix/var}$}{Fix/variable cost of a unit of generator $g$.}
\nomenclature[P]{${c}_g^\text{su/sd}$}{Startup/shutdown cost of a unit of generator $g$.}
\nomenclature[P]{${p}_{c,t}^{\text{}}$}{Load demand of consumer $c$.}
\nomenclature[P]{${p}_{p,t}^{\text{}}$}{Load demand of prosumer $p$.}
\nomenclature[P]{${p}_{n,t}^{\text{r}}$}{Power reserve requirement of node $n$.}
\nomenclature[P]{$\underline{{x}}/\overline{{x}}$}{Minimum/maximum limit of variable $x$.}
\nomenclature[P]{$\overline{{U}}_{g}$}{Total number of identical units of generator $g$.}
\nomenclature[P]{${{r}}^{+/-}_{g}$}{Ramp-up/down rate of a unit of generator $g$.}
\nomenclature[P]{$\tau^{\text{u/d}}_{g}$}{Minimum up/down time of a unit of generator $g$.}
\nomenclature[P]{$\tilde{t}$}{Time slot offset index.}
\nomenclature[P]{$\Delta{{t}}$}{Time resolution.}
\nomenclature[P]{${B}_l$}{Susceptance of line $l$.}
\nomenclature[P]{${p}_{g,t}^{\text{RES}}$}{Max. output power of renewable generator $g \in \mathcal{G}^\text{RES}$.}
\nomenclature[P]{${p}_{g,t}^{\text{CST}}$}{Max. thermal power capture by generator $g \in \mathcal{G}^\text{CST}$.}
\nomenclature[P]{${H}_g$}{Inertia of a unit of generator $g$.}
\nomenclature[P]{${S}_g$}{MVA rating of a unit of generator $g$.}
\nomenclature[P]{${H}_{n,t}$}{Minimum synchronous inertia requirement of node $n$.}
\nomenclature[P]{$\eta_x$}{Efficiency of component $x$.}
\nomenclature[P]{${p}_{p,t}^{\text{pv}}$}{Aggregated PV power of prosumer $p$.}
\nomenclature[P]{$\lambda$}{Feed-in price ratio.}

\nomenclature[I]{$\hat{{s}}_{g}$}{Number of online units of generator $g \in \mathcal{G}^\text{syn}$ at start of horizon.}
\nomenclature[I]{$\hat{{p}}_{g}$}{Power dispatch of generator $g$ at start of horizon.}
\nomenclature[I]{$\hat{{u}}_{g,t}$}{Minimum number of units of generator $g \in \mathcal{G}^\text{syn}$ required to remain online for time $t<\tau_g^\text{u}$.}
\nomenclature[I]{$\hat{{d}}_{g,t}$}{Minimum number of units of generator $g \in \mathcal{G}^\text{syn}$ required to remain offline for time $t<\tau_g^\text{d}$.}
\nomenclature[I]{$\hat{{e}}_{g}$}{Energy stored in TES of $g \in \mathcal{G^\text{CST}}$ at start of horizon.}
\nomenclature[I]{$\hat{{e}}_{s}$}{Energy stored in storage plant $s$ at start of horizon.}
\nomenclature[I]{$\hat{{e}}_{p}^\text{b}$}{Battery state of charge for prosumer $p$ at start of horizon.}

%\makenomenclature
\printnomenclature

\IEEEpeerreviewmaketitle

%\vspace{-1.0em}
\section{Introduction}
Power systems worldwide are moving away from domination by large-scale synchronous generation and passive consumers.
Instead, in future grids\footnote{We interpret a \textit{future grid} to mean the study of national grid type structures with the transformational changes over the long-term out to 2050.} new actors, such as variable renewable energy sources (RES)\footnote{For the sake of brevity, by RES we mean ``unconventional'' renewables like wind and solar, but excluding conventional RES, like hydro, and dispatchable unconventional renewables, like concentrated solar thermal.}, price-responsive users equipped with small-scale PV-battery systems (called \emph{prosumers}), demand response (DR), and energy storage will play an increasingly important role.
Given this, in order for policy makers and power system planners to evaluate the integration of high-penetrations of these new elements into future grids, new simulation tools need to be developed. 
Specifically, there is a pressing need to understand the effects of technological change on future grids, in terms of energy balance, stability, security and reliability, over a wide range of highly-uncertain future scenarios.
This is complicated by the inherent and unavoidable uncertainty surrounding the availability, quality and cost of new technologies (e.g. battery or photo-voltaic system costs, or concentrated solar thermal (CST) generation operating characteristics) and the policy choices driving their uptake.
The recent blackout in South Australia \cite{AEMO2016} serves as a reminder that things can go wrong when the uptake of new technologies is not planned carefully. 

Future grid planning thus requires a major departure from conventional power system planning, where only a handful of the most critical scenarios are analyzed. To account for a wide range of possible future evolutions, \emph{scenario analysis} has been proposed in many industries, e.g. in finance and economics \cite{LearningFromTheFuture_1998}, and in energy \cite{Foster2013, Sanchis2015}. In contradistinction to power system planning, where the aim is to find an optimal transmission and/or generation expansion plan, the aim of scenario analysis is to analyze possible evolution pathways to inform power system planning and policy making. Given the uncertainty associated with long-term projections, the focus of future grid scenario analysis is limited only to the analysis of what is technically possible, although it might also consider an explicit costing \cite{Elliston2016}.
In more detail, existing future grid feasibility studies have shown that the balance between demand and supply can be maintained even with high penetration of RESs by using large-scale storage, flexible generation, and diverse RES technologies \cite{Energy2010,  Elliston2013, Budischak2013, Mason2010, AEMORES}.
However, they only focus on balancing and use simplified transmission network models (either copper plate or network flow; a notable exception is the Greenpeace pan-European study \cite{GP} that uses a DC load flow model). This ignores network related issues, which limits these models' applicability for stability assessment. 

\begin{figure}
	\centering
		\includegraphics[width=85mm] {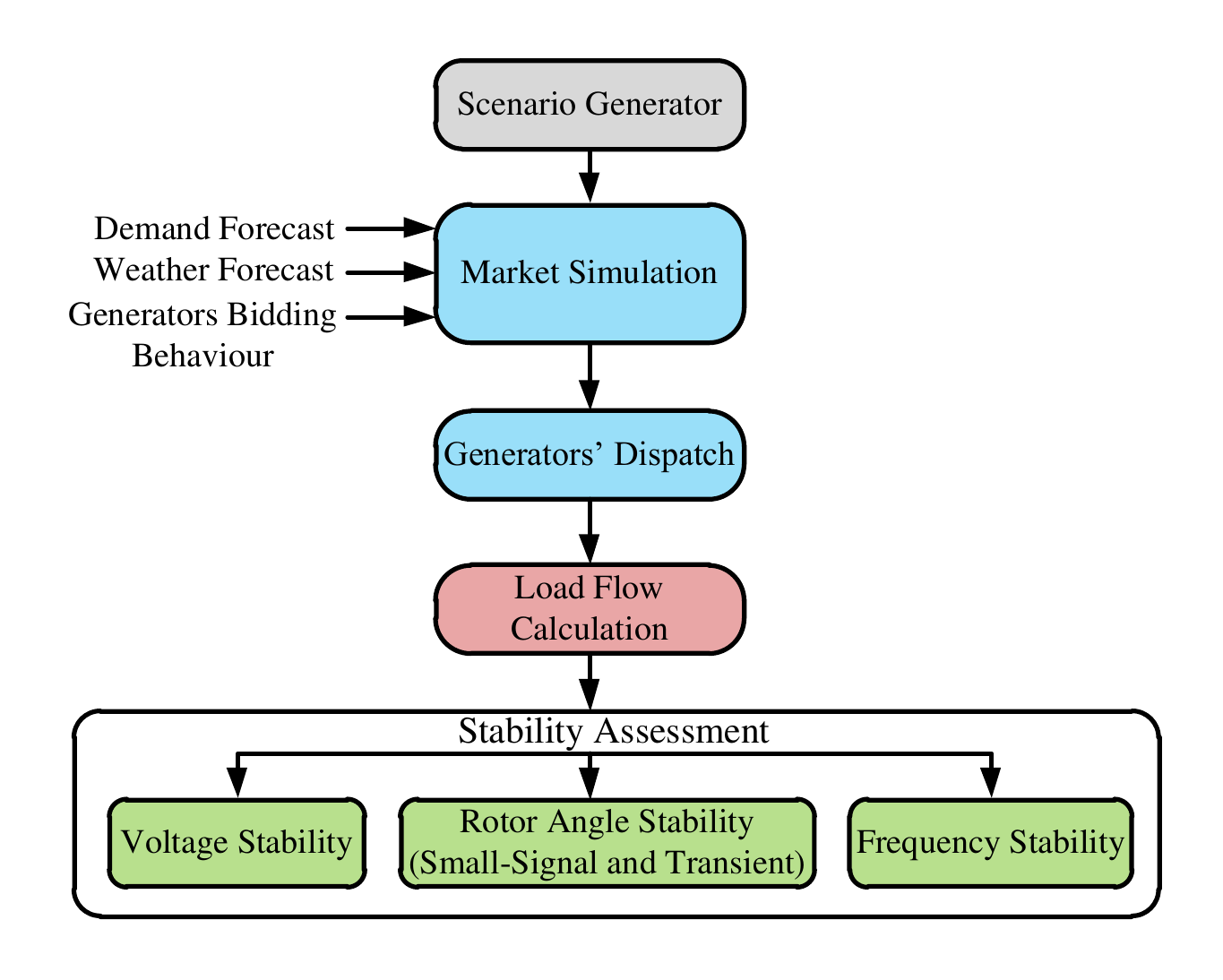}
		%\vspace{-1.0em}
	\caption{Simulation platform for the performance and stability assessment of future grid scenarios.}
	\label{fig:future grid_fig}
	%\vspace{-1.5em}
\end{figure}

To the best of our knowledge, the Future Grid Research Program, funded by the Australian Commonwealth Scientific and Industrial Research Organisation (CSIRO) is the first to propose a comprehensive modeling framework for future grid scenario analysis that also includes \emph{stability assessment}. The aim of the project is to explore possible future pathways for the evolution of the Australian grid out to 2050 by looking beyond simple balancing. To this end, a simulation platform has been proposed in \cite{Marzooghi2014} that consists of a market model, power flow analysis, and stability assessment, Fig.~\ref{fig:future grid_fig}. The platform has been used, with additional improvements, to study fast stability scanning \cite{Liu2016}, inertia \cite{Ahmad16a}, modeling of prosumers for market simulation \cite{Marzooghi16b, Marzooghi16},  impact of prosumers on voltage stability \cite{Riaz16b}, and power system flexibility using CST \cite{Riaz16a} and battery storage \cite{Riaz15}.
In order to capture the inter-seasonal variations in the renewable generation, computationally intensive time-series analysis needs to be used.
A major computational bottleneck of the framework is the market simulation.
%The proposed electricity market simulation tool (MST) is based on a unit commitment (UC) problem that is used to generate generator dispatch traces. UC has always played a vital role in power systems operation and planning. In an operational framework, reliability UC is used for ensuring security and day-ahead market clearance. For planning studies, UC provides a powerful and flexible tool for investigating the medium-term impacts of investment decisions. To reduce the computational complexity, several approximations are used, depending on the application.

Within this context, the contribution of this paper is to propose a unified generic market simulation tool (MST) based on a unit commitment (UC) problem suitable for future grid scenario analysis, including stability assessment. The tool incorporates the following key features:
\begin{itemize}
	\item market structure agnostic modeling framework,
	\item integration of various types and penetrations of RES and emerging demand-side technologies,
	\item generic demand model considering the impact of prosumers,
	\item explicit network representation, including HVDC lines, using a DC power flow model, 
	\item explicit representation of the number of online synchronous generators,
	\item explicit representation of system inertia and reactive power support capability of synchronous generators,
	\item computational efficiency with sufficient accuracy.
\end{itemize}

The presented model builds on our existing research~\cite{Marzooghi16,Marzooghi16b,Riaz15,Riaz16b,Ahmad16a,Riaz16a} and combines all these in a single coherent formulation.

In more detail, to reduce the computational burden, the following techniques are used building on the methods proposed in~\cite{Carrion06,Palmintier14a}:
\begin{itemize}
	\item unit clustering,
	\item rolling horizon approach,
	\item constraint clipping.
\end{itemize} 

The computational advantages of our proposed model are shown on a simplified 14-generator model of the Australian National Energy Market (NEM) as a test grid~\cite{Gibbard10}. Four cases for different RES penetration are run for one to seven days horizon length, and computational metrics are reported. To reflect the accuracy of the proposed MST, system inertia and voltage stability margins are used as a benchmark. In simulations, RES and load traces are taken from the National Transmission Network Developed Plan (NTNDP) data, provided by the Australian Energy Market Operator (AEMO)~\cite{AEMO2016b}.

The remainder of the paper is organized as follows: Literature review and related work are discussed in Section II, while Section
III details the MST. A detailed description
of the simulation setup is given in Section IV. In Section V results are analyzed and discussed in detail. Finally, Section VI concludes the paper.

%\vspace{-1.0em}
\section{Related Work}
In order to better explain the functional requirements of the proposed MST, we first describe the canonical UC formulation.
An interested reader can find a comprehensive literature survey in \cite{Tahanan2015}.

%\vspace{-1.0em}
\subsection{Canonical Unit Commitment Formulation}
The UC problem is an umbrella term for a large class of problems in power system operation and planning whose objective is to schedule and dispatch power generation at minimum cost to meet the anticipated demand, while meeting a set of system-wide constraints. In smart grids, problems with a similar structure arise in the area of energy management, and they are sometimes also called UC \cite{Ramos14}.
Before deregulation, UC was used in vertically integrated utilities for generation scheduling to minimize production costs. After deregulation, UC has been used by system operators to maximize social welfare, but the underlying optimization model is essentially the same.

Mathematically, UC is a large-scale, nonlinear, mixed-integer optimization problem under uncertainty. With some abuse of notation, the UC optimization problem can be represented in the following compact formulation \cite{Bertsimas2013}:
\begin{align}
\mathop{\operatorname{minimize}}\limits_{ \mathbf{x}_{\text{c}}, \mathbf{x}_{\text{b}} } \quad & f_{\text{c}}(\mathbf{x}_{\text{c}}) + f_{\text{b}}(\mathbf{x}_{\text{b}}) \label{obj} \\
\mathop{\operatorname{subject\; to}} \quad & g_{\text{c}}(\mathbf{x}_{\text{c}}) \le \mathbf{b} \label{const_real} \\ 
& g_{\text{b}}(\mathbf{x}_{\text{b}}) \le \mathbf{c} \label{const_bin} \\
& h_{\text{c}}(\mathbf{x}_{\text{c}}) + h_{\text{b}}(\mathbf{x}_{\text{b}})\le \mathbf{d} \label{const_realbin} \\
& \mathbf{x}_{\text{c}} \in \mathbb{R}^{+},\: \mathbf{x}_{\text{b}} \in \{ 0,1 \} \nonumber
\end{align}

Due to the time-couplings, the UC problem needs to be solved over a sufficiently long horizon.
The decision vector $\mathbf{x} = \{ \mathbf{x}_{\text{c}}, \mathbf{x}_{\text{b}} \}$ for each time interval consist of continuous and binary variables. The continuous variables, $\mathbf{x}_{\text{c}}$, include  generation dispatch levels, load levels, transmission power flows, storage levels, and transmission voltage magnitudes and phase angles. The binary variables, $\mathbf{x}_{\text{b}}$, includes scheduling decisions for generation and storage, and logical decisions that ensure consistency of the solution.
The objective (\ref{obj}) captures the total production cost, including fuel costs, start-up costs and shut-down costs.
The constraints include, respectively: dispatch related constraints such as energy balance, reserve requirements, transmission limits, and ramping constraints (\ref{const_real}); commitment variables, including minimum up and down, and start-up/shut-down constraints (\ref{const_bin}); and constraints coupling commitment and dispatch decisions, including minimum and maximum generation capacity constraints (\ref{const_realbin}).  

The complexity of the problem stems from the following: (i) certain generation technologies (e.g. coal-fired steam units) require long start-up and shut-down times, which requires a sufficiently long solution horizon; (ii) generators are interconnected, which introduces couplings through the power flow constraints; (iii) on/off decisions introduce a combinatorial structure; (iv) some constraints (e.g. AC load flow constraints) and parameters (e.g. production costs) are non-convex; and (v) the increasing penetration of variable renewable generation and the emergence of demand-side technologies introduce uncertainty.
As a result, a complete UC formulation is computationally intractable, so many approximations and heuristics have been proposed to strike a balance between computational complexity and functional requirements. For example, power flow constraints can be neglected altogether (a copper plate model), can be replaced with simple network flow constraints to represent critical inter-connectors, or, instead of (non-convex) AC, a simplified (linear) DC load flow is used. 

%\vspace{-1.0em}
\subsection{UC Formulations in Existing Future Grid Studies}
In operational studies: the nonlinear constraints, e.g. ramping, minimum up/down time (MUDT) and thermal limits are typically linearized; startup and shutdown exponential costs are discretized, and; non-convex and non-differentiable variable cost functions are expressed as piecewise linear function \cite{Arroyo00,Carrion06}. In planning studies, due to long horizon lengths, the UC model is simplified even further. For example: combinatorial structure is avoided by aggregating all the units installed at one location \cite{Palmintier14a,Hara1966,Langrene2011}; piecewise linear cost functions and constraints are represented by one segment only; some costs (e.g. startup, shutdown and fix costs) are ignored; a deterministic UC with perfect foresight is used, and; non-critical binding constraints are omitted~\cite{Tuohy09,Palmintier11}\footnote{An interested reader can refer to \cite{Palmintier14b} for a discussion on binding constraints elimination for generation planning.}.% \textcolor{blue}{or replaced with simple restrictions}. 
To avoid the computational complexity associated with the mixed integer formulation, a recent work~\cite{Zhang16} has proposed a linear relaxation of the UC formulation for flexibility studies, with an accuracy comparable to the full binary mixed integer linear formulation. 

In contrast to operation and planning studies, the computational burden of future grid scenario analysis is even bigger, due to a sheer number of scenarios that need to be analyzed, which requires further simplifications. For example, the Greenpeace study \cite{GP} uses an optimal power flow for generation dispatch and thus ignores UC decisions. Unlike the Greenpeace study, the Irish All Island Grid Study~\cite{IRE} and the European project e-Highway2050 \cite{Ehighway2050} ignore load flow constraints altogether, however they do use a rolling horizon UC, with simplifications. The Irish study, for example doesn't put any restriction on the minimum number of online synchronous generators to avoid RES spillage, and the e-Highway2050 study uses a heuristics to include DR. The authors of the e-Highway2050 study, however, acknowledge the size and the complexity of the optimization framework in long term planning, and plan to develop new tools with a simplified network representation~\cite{Ehighway2050}.

In summary, a UC formulation depends on the scope of the study. Future grid studies that explicitly include stability assessment bring about some specific requirements that are routinely neglected in the existing UC formulations, as discussed next. 

%%\vspace{-1.0em}
\section{Market Simulation Tool}
\subsection{Functional Requirements}
The focus of our work is stability assessment of future grid scenarios. Thus, MST must produce dispatch decisions that accurately capture the kinetic energy stored in rotating masses (inertia), active power reserves and reactive power support capability of synchronous generators, which all depend upon the number of online units and the respective dispatch levels. 

\begin{figure}
	\centering
	\includegraphics[width=71mm] {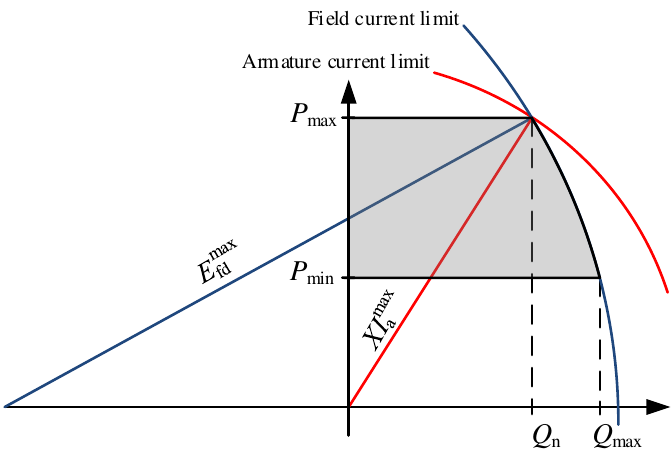}
	%\vspace{-1.0em}
	\caption{Illustrative operating chart a synchronous generator in an over-excited mode (shaded region).}
	\label{fig:PQcurve}
%	%\vspace{-1.5em}
\end{figure}

\begin{figure}
	\centering
	\includegraphics[width=85mm] {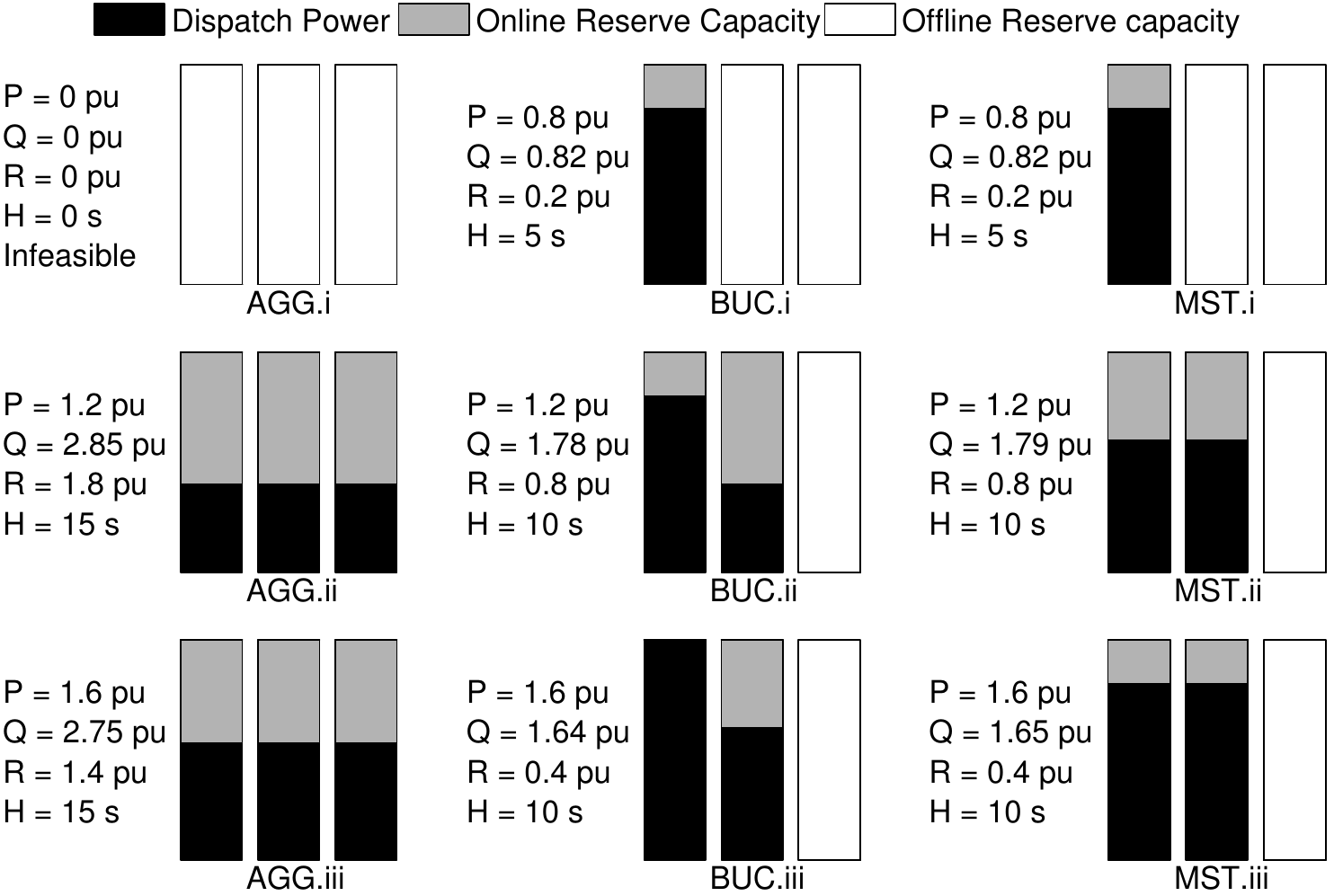}
	%\vspace{-1.0em}
	\caption{Illustrative operating cases, (i) $P = \SI{0.8}{pu}$, (ii) $P = \SI{1.2}{pu}$, (iii) $P = \SI{1.6}{pu}$, showing active power dispatch level $P$, reactive power support capability $Q$, online active power reserves $R$ and generator inertia $H$. The columns correspond to three different UC formulations: aggregated (AGG), binary (BUC) and the proposed MST formulation.}
	\label{fig:Generation_comparision}
%	%\vspace{-1em}
\end{figure}

For the sake of illustration, consider a generation plant consisting of three identical (synchronous) thermal units, with the following characteristics: (i) constant terminal voltage of \SI{1}{pu}; (ii) minimum technical limit $P_{\textrm{min}} = \SI{0.4}{pu}$; (iii) power factor of $0.8$; (iv) maximum excitation limit $E_{\textrm{fd}}^{\textrm{max}} = \SI{1.5}{pu}$; and (v) normalized inertia constant $H = \SI{5}{\second}$. We further assume that in the over-excited region, the excitation limit is the binding constraint, as shown in Fig. \ref{fig:PQcurve}. Observe that the maximum reactive power capability depends on the active power generated, and varies between $Q_{\textrm{n}}$ at $P_{\textrm{max}} = \SI{1}{pu}$ and $Q_{\textrm{max}}$ at $P_{\textrm{min}}$.
We consider three cases defined by the total active power generation of the plant: (i) \SI{0.8}{pu}, (ii) \SI{1.2}{pu}, and (iii) \SI{1.6}{pu}. 
The three scenarios correspond to the rows  in Fig.~\ref{fig:Generation_comparision}, which shows the active power dispatch level $P$, reactive power support capability $Q$, online active power reserves $R$, and generator inertia $H$. 
The three columns show feasible solutions for three different UC formulations: all three units are aggregated into one equivalent unit (AGG), standard binary UC (BUC) when each unit is modeled individually, and the proposed market simulation tool (MST). A detailed comparison of the three formulations is given in Section V.

Although the results are self-explanatory, a few things are worth emphasizing. In case (i), aggregating the units into one equivalent unit (AGG) results in the unit being shut down due to the minimum technical limit. The individual unit representation (BUC), on the other hand, does allow the dispatch of one or two units, but with significantly different operational characteristics. In cases (ii) and (iii), the total inertia in the AGG formulation is much higher, which has important implications for frequency stability. A similar observation can be made for the reactive power support capability, which affects voltage stability. Also, dispatching power from all three units results in a significantly higher active power reserve. And last, a higher reactive power generation due to a lower $P$ reduces the internal machine angle, which improves transient stability.

In conclusion, a faithful representation of the number of online synchronous machines is of vital importance for stability assessment. An individual unit representation, however, is computationally expensive, so the computational burden should be reduced, as discussed in the following section. Next, an explicit network representation is required. An AC load flow formulation, however, is nonlinear (and non-convex), which results in an intractable mixed-integer nonlinear problem. Therefore, we use a DC load flow representation with a sufficiently small voltage angle difference on transmission lines. Our experience shows that an angle difference of \SI{30}{\degree} results in a manageable small number of infeasible operating conditions that can be dealt with separately.

%\vspace{-1.0em}
\subsection{Computational Speedup} \label{Section:RT_Improve}
The MST is based on the UC formulation using constant fixed, startup, shutdown and production costs. To improve its computational efficiency, the dimensionality of the optimization problem is reduced employing: (i) unit clustering \cite{Palmintier14a} to reduce the number of variables needed to represent a multi-unit generation plant; (ii) a rolling horizon approach \cite{Tuohy09,Ramos14,IRE-2B-2008} to reduce the time dimension; and (iii) constraint clipping to remove most non-binding constraints.  

\subsubsection{Unit Clustering}
Linearized UC models are computationally efficient for horizons of up to a few days, which makes them extremely useful for operational studies. For planning studies, however, where horizon lengths can be up to a year, or more, these models are still computationally too expensive. Our work builds on the clustering approach proposed in~\cite{Palmintier14a}, where identical units at each generation plant are aggregated by replacing binary variables with fewer integer variables. The status of online units, startup/shutdown decisions and dispatched power are tracked by three integer variables and one continuous variable \emph{per plant} per period, as opposed to three binary and one continuous variable \emph{per unit} per period. Further clustering proposed in~\cite{Palmintier14a} is not possible in our formulation because of the explicit network representation required in the MST. %Henceforth, the UC model with binary decision variables will be referred to as binary UC (BUC) and the UC with the clustering approach will be referred to as integer UC (IUC). 

\subsubsection{Rolling Horizon} 
Solving the UC as one block, especially for long horizons, is computationally too expensive. This can be overcome by breaking the problem into several smaller intervals called sub-horizons~\cite{Tuohy09,Ramos14,IRE-2B-2008}. To ensure accuracy and consistency of the solution, a proper overlap between sub-horizons is maintained and the terminating state of the previous sub-horizon is used as the initial condition of the next sub-horizon. The minimum sub-horizon length depends on the time constants associated with the decision variables. While these might be in the order of hours for thermal power plants, they can be significantly longer for energy storage. Large-scale hydro dams, for example, require horizon lengths of several weeks, or even months. In our research, however, the sub-horizon length is up to a few days to cater for thermal energy storage (TES) of CST plants and battery storage. The optimization of hydro dams is not explicitly considered, however it can be taken into account heuristically, if needed.

\subsubsection{Constraint Clipping}
The size of the problem can be reduced by removing non-binding constraints, which doesn't affect the feasible region. For instance, an MUDT constraint on a unit with an MUDT less than the time interval is redundant\footnote{This is especially the case when the time resolution is coarse. In our studies, the time step is one hour. In operational studies, where the resolution can be as short as five minutes, constraint clipping is less useful.}. Similarly, a ramp constraint for flexible units is redundant if the time step is sufficiently long. With a higher RES penetration, in particular, where backup generation is provided by fast-ramping gas turbines, this technique can significantly reduce the size of the optimization problem, and hence improves the computational performance due to a larger number of units with higher ramp rates and smaller MUDTs. It should be noted that optimization pre-solvers might not able to automatically remove these constraints. 

%%%\vspace{-1.0em}
\subsection{MST UC Formulation}

\subsubsection{Objective function}
The objective of the proposed MST is to minimize total generation cost for all sub-horizons $h$:
\begin{multline}
\underset{\Omega}{\mathop{\operatorname{minimize}}}  \  \sum_{t\in \mathcal{T}}^{} \sum_{g\in \mathcal{G}}^{}  
\Big( {c}_g^\text{fix} s_{g,t}
+{c}_g^\text{su} u_{g,t} + {c}_g^\text{sd}  d_{g,t} +{c}_g^\text{var}  p_{g,t} \Big), 
\end{multline}	
where $\Omega = \{s_{g,t},u_{g,t},d_{g,t},p_{g,t}, p_{s,t}, p_{l,t} \}$ are the decision variables of the problem, and ${c}_g^\text{fix}$, ${c}_g^\text{su}$, ${c}_g^\text{sd}$, and ${c}_g^\text{var}$ are fixed, startup, shutdown and variable cost, respectively.
As typically done in planning studies \cite{Palmintier14a}, \cite{Zhang16}, the costs are assumed constant to reduce the computational complexity. The framework, however, also admits a piece-wise linear approximation proposed in \cite{Carrion06}.

\subsubsection{System Constraints}
System Constraints\footnote{All the constraints must be satisfied in all time slots $t$, however, for sake of notational brevity, this is not explicitly mentioned.} include power balance constraints, power reserve and minimum synchronous inertia requirements.

\textbf{Power balance:}
Power generated at node $n$ must be equal to the node power demand plus the net power flow on transmission lines connected to the node:
\begin{multline}
\sum_{g\in \mathcal{G}_n}^{}{p_{g,t}} =  \\
\sum_{c \in \mathcal{C}_n}p_{{c},t}^{} + \sum_{p \in \mathcal{P}_n} p_{{p},t}^{\text{g}+} - \sum_{p \in \mathcal{P}_n} p_{{p},t}^{\text{g}-} + \sum_{s \in \mathcal{S}_n}{p}_{s,t} + \sum_{l \in \mathcal{L}_n}(p_{l,t} +\Delta p_{l,t}), 
\end{multline}
where $ \mathcal{G}_n, \mathcal{C}_n, \mathcal{P}_n, \mathcal{S}_n, \mathcal{L}_n$ represent respectively the set of generators, consumers, prosumers\footnote{Price-responsive users equipped with small-scale PV-battery systems.}, utility storage plants and lines connected to node $n$. 

\textbf{Power reserves:}
To cater for uncertainties, active power reserves provided by synchronous generation $g \in \mathcal{G}^\text{syn}$ are maintained in each region $r$:
\begin{multline}
\sum_{g \in \{ (\mathcal{G}^\text{syn}-\mathcal{G}^\text{CST}) \cap \mathcal{G}^r\} } (\overline{{p}}_{g} s_{g,t} - {p}_{g,t}) + \\\sum_{g \in \{ \mathcal{G}^\text{CST} \cap \mathcal{G}^r\}} \mathop{\operatorname{min}}(\overline{{p}}_{g} s_{g,t} - {p}_{g,t},{e}_{g,t}-{p}_{g,t})  \geq \sum_{n \in \mathcal{N}_r}^{}{p}_{n,t}^{\text{r}}. \label{eq:3}
\end{multline}
For synchronous generators other than concentrated solar thermal (CST), reserves are defined as the difference between the online capacity and the current operating point. For CST, reserves can either be limited by their online capacity or energy level of their thermal energy system (TES).
Variable $s_{g,t}$ in (\ref{eq:3}) represents the total number of online units at each generation plant, and $\mathcal{G}^r$ and $\mathcal{N}_r$ represent the sets of generators and nodes in region $r$, respectively.

\textbf{Minimum synchronous inertia requirement:} 
To ensure frequency stability, a minimum level of inertia provided by synchronous generation must be maintained at all times (more details are available in \cite{Ahmad16a}) in each region $r$:  
\begin{equation}
\sum_{g \in \{\mathcal{G}^\text{syn} \cap \mathcal{G}^r\} }^{} s_{g,t} {H}_g {S}_g \geq \sum_{n \in \mathcal{N}_r}^{}{H}_{n,t}.
\end{equation}

\subsubsection{Network constraints}
Network constraints include DC power flow constraints and thermal line limits for AC lines, and active power limits for HVDC lines.  

\textbf{Line power constraints:} 
A DC load flow model is used for computational simplicity for AC transmission lines\footnote{A sufficiently small ($\sim\!\SI{30}{\degree}$) voltage angle difference over a transmission line is used to reduce the number of nonconvergent AC power flow cases.}:
\begin{equation}
p_{l,t}^{x,y} = {B}_l(\delta_{x,t} - \delta_{y,t}), \quad l \in  \mathcal{L^\text{AC}},
\end{equation}
where the variables $\delta_{x,t}$ and $\delta_{y,t}$ represent voltage angles at nodes $x \in \mathcal{N}$ and $y \in \mathcal{N}$, respectively.

\textbf{Thermal line limits:}
Power flows on all transmission lines are limited by the respective thermal limits of line $l$:
\begin{equation}
\mid  p_{l,t} \mid \leq \overline{{p}}_l,
\end{equation}
where $\overline{{p}}_l$ represents the thermal limit of line $l$.

\subsubsection{Generation constraints}
Generation constraints include physical limits of individual generation units.
For the binary unit commitment (BUC), we adopted a UC formulation requiring three binary variables per time slot (on/off status, startup, shutdown) to model an individual unit. In the MST, identical units of a plant are clustered into one individual unit~\cite{Palmintier14a}. 
This requires three \emph{integer} variables (on/of status, startup, and shutdown) \emph{per generation plant} per time slot as opposed to three \emph{binary} variables \emph{per generation unit} per time slot in the BUC, as discussed in Section III.B of \emph{A Computationally Efficient Market Model for Future Grid Scenario Studies.}

\textbf{Generation limits:} 
Dispatch levels of a synchronous generator $g$ are limited by the respective stable operating limits:
\begin{equation}
s_{g,t} \underline{{p}}_{g} \leq p_{g,t} \leq s_{g,t} \overline{{p}}_{g}, \quad  g \in \mathcal{G}^\text{syn}.
\end{equation}
The power of RES\footnote{For the sake of brevity, by RES we mean ``unconventional'' renewables like wind and solar, but excluding conventional RES, like hydro, and dispatchable unconventional renewables, like concentrated solar thermal.} generation is limited by the availability of the corresponding renewable resource (wind or sun):
\begin{equation}
s_{g,t} \underline{{p}}_{g} \leq p_{g,t} \leq s_{g,t} {{p}}_{g,t}^{\text{RES}}, \quad  g \in \mathcal{G}^\text{RES}.
\end{equation}

\textbf{Unit on/off constraints:} 
A unit can only be turned on if and only if it is in off state and vice versa:
\begin{equation}
u_{g,t}-d_{g,t}=s_{g,t}-s_{g,t-1}, \quad  t \neq 1, \ g \in  \mathcal{G}^\text{syn}. \label{eq:9}
\end{equation}	
In a rolling horizon approach, consistency between adjacent time slots is ensured by:
\begin{equation}
u_{g,t}-d_{g,t}=s_{g,t} - \hat{{s}}_{g}, \quad  t =1, \ g \in  \mathcal{G}^\text{syn}, \label{eq:10}
\end{equation}
where $\hat{{s}}_{g}$ is the initial number of online units of generator $g$. Equations (\ref{eq:9}) and (\ref{eq:10}) also implicitly determine the upper bound of $u_{g,t}$ and $d_{g,t}$ in terms of changes in  $s_{g,t}$.

\textbf{Number of online units:}
Unlike the BUC, the MST requires an explicit upper bound on status variables:
\begin{equation}
s_{g,t} \leq \overline{ {U}}_{g}. 
\end{equation}

\textbf{Ramp-up and ramp-down limits:}
Ramp rates of synchronous generation should be kept within the respective ramp-up (\ref{eg:12}), (\ref{eg:13}) and ramp-down limits (\ref{eg:14}), (\ref{eg:15}):
\begin{align}
& p_{g,t} - p_{g,t-1} \leq s_{g,t} {r}^+_g, & t \neq 1, g \in  \{\mathcal{G}^\text{syn} | {r}^+_g < \overline{ {p}}_{g}\}, \label{eg:12} \\ 
& p_{g,t} - \hat{p}_{g} \leq s_{g,t} {r}^+_g, & t =1, g \in  \{\mathcal{G}^\text{syn} | {r}^+_g < \overline{ {p}}_{g}\},  \label{eg:13} \\			
& \hat{p}_{g} - p_{g,t} \leq s_{g,t-1} {r}^-_g,  & t \neq 1, g \in  \{\mathcal{G}^\text{syn} | {r}^-_g < \overline{ {p}}_{g}\},  \label{eg:14} \\
& \hat{p}_{g} - p_{g,t} \leq \hat{{s}}_{g} {r}^-_g, & t =1, g \in  \{\mathcal{G}^\text{syn} | {r}^-_g < \overline{ {p}}_{g}\}.  \label{eg:15} 
\end{align}
In the MST, a ramp limit of a power plant is defined as a product of the ramp limit of an individual unit and the number of online units in a power plant $s_{g,t}$. If $s_{g,t}$ is binary, these ramp constraints are mathematically identical to ramp constraints of the BUC.   
If a ramp rate multiplied by the length of the time resolution $\Delta{{t}}$ is less than the rated power, the rate limit has no effect on the dispatch, so the corresponding constraint can be eliminated.
Constraints explicitly defined for $t=1$ are used to join two adjacent sub-horizons in the rolling-horizon approach.

\textbf{Minimum up and down times:} 
Steam generators must remain on for a period of time  $\tau_{g}^\text{u}$ once turned on (minimum up time):
\begin{align}
& s_{g,t} \geq \sum_{\tilde{t}=\tau_{g}^\text{u}-1}^{0} u_{g,t-\tilde{t}}, & t \geq \tau_{g}^\text{u}, \  g \in  \{\mathcal{G}^\text{syn} | \tau_{g}^\text{u} > \Delta{{t}}\}, \\
& s_{g,t} \geq \sum_{\tilde{t}=t-1}^{0} {u_{g,t-\tilde{t}}}  + \hat{{u}}_{g,t}, & t < \tau_{g}^\text{u}, \ g \in  \{\mathcal{G}^\text{syn} | \tau_{g}^\text{u} > \Delta{{t}}\}.
\end{align}		
Similarly, they must not be turned on for a period of time  $\tau_{g}^\text{d}$ once turned off (minimum down time):
\begin{align}
& s_{g,t} \leq \overline {{U}}_{g} - \sum_{\tilde{t}=\tau_{g}^\text{d}-1}^{0} d_{g,t-\tilde{t}}, & t \geq \tau_{g}^\text{d}, \ g \in  \{\mathcal{G}^\text{syn} | \tau_{g}^\text{d} > \Delta{{t}}\}, \\
& s_{g,t} \leq \overline {{U}}_{g} - \sum_{\tilde{t}=t-1}^{0}{d_{g,t-\tilde{t}} } - \hat{{d}}_{g,t}, & t < \tau_{g}^\text{d}, \ g \in  \{\mathcal{G}^\text{syn} | \tau_{g}^\text{d} > \Delta{{t}}\}.
\end{align}	
Similar to the rate limits, if the minimum up and down times are smaller than the time resolution $\Delta{{t}}$, the corresponding constraints can be eliminated. 
Due to integer nature of discrete variables in the MST, the definition of the MUDT constraints in the RH approach requires the number of online units for the last $\tau^{\text{u/d}}$ time interval to establish the relationship between the adjacent sub-horizons. If the $\tau_{g}^{\text{u/d}}$ is smaller than time resolution $\Delta{{t}}$, then these constraints can be eliminated. 

\subsubsection{CST constraints:}
CST constraints include TES energy balance and storage limits.                

\textbf{TES state of charge (SOC)} 
determines the TES energy balance subject to the accumulated energy in the previous time slot, thermal losses, thermal power provided by the solar farm and electrical power dispatched from the CST plant:
\begin{align}
& e_{g,t}=\eta_{g}{e}_{g,t-1}+{{p}}_{g,t}^{\text{CST}}-p_{g,t},	& t \neq 1, \ g \in  \mathcal{G}^\text{CST}, \\
& e_{g,t}=\eta_{g}\hat{{e}}_{g}+{{p}}_{g,t}^{\text{CST}}-p_{g,t},  & t=1, \ g \in  \mathcal{G}^\text{CST},
\end{align}	
where, ${{p}}_{g,t}^{\text{CST}}$ is the thermal power collected by the solar field of generator $g \in \mathcal{G}^\text{CST}$. 

\textbf{TES limits:}	Energy stored is limited by the capacity of a storage tank:
\begin{equation}
\underline{e}_{g}\leq{e}_{g,t} \leq \overline{e}_{g},	\quad g \in  \mathcal{G}^\text{CST}.
\end{equation}

\subsubsection{Utility storage constraints}
Utility-scale storage constraints include energy balance, storage capacity limits and power flow constraints. The formulation is generic and can capture a wide range of storage technologies.  

\textbf{Utility storage SOC limits} determine the energy balance of storage plant  $s$:
\begin{align}
& e_{s,t}=\eta_{s}{e}_{s,t-1}+{p}_{s,t},	& t \neq 1, \\
& e_{s,t}=\eta_{s}\hat{{e}}_{s}+{p}_{s,t},	& t=1.
\end{align}	

\textbf{Utility storage capacity limits:}
Energy stored is limited by the capacity of storage plant $s$:
\begin{equation}
\underline{e}_{s}\leq{e}_{s,t} \leq \overline{e}_{s}.
\end{equation}	

\textbf{Charge/discharge rates} limit the charge and discharge powers of storage plant $s$:
\begin{equation}
\overline{{p}}_{s}^- \leq {p}_{s,t} \leq \overline{{p}}_{s}^+,
\end{equation}
where $\overline{{p}}_{s}^-$ and $\overline{{p}}_{s}^+$ represent the maximum power discharge and charge rates of a storage plant, respectively.

\subsubsection{Prosumer sub-problem}
The prosumer sub-problem captures the aggregated effect of prosumers. It is modeled using a bi-level framework in which the upper-level unit commitment problem described above minimizes the total generation cost, and the lower-level problem maximizes prosumers' self-consumption. The coupling is through the prosumers' demand, not through the electricity price, which renders the proposed model market structure agnostic. As such, it implicitly assumes a mechanism for demand response aggregation. The Karush-Kuhn-Tucker optimality conditions of the lower-level problem are added as the constraints to the upper-level problem, which reduces the problem to a single mixed integer linear program.

The model makes the following assumptions: (i) the loads are modeled as price anticipators; (ii) the demand model representing an aggregator consists of a large population of prosumers connected to an unconstrained distribution network who collectively maximize self-consumption; (iii) aggregators do not alter the underlying power consumption of the prosumers; and (iv) prosumers have smart meters equipped with home energy management systems for scheduling of the PV-battery systems, and, a communication infrastructure is assumed that allows a two-way communication between the grid, the aggregator and the prosumers. More details can be found in \cite{Marzooghi16b}.

\textbf{Prosumer Objective function:} 
Prosumers aim to minimize electricity expenditure:
\begin{equation}
\underset{p_{p}^{g\text{+/--}}, p_{p}^{b}}{\mathop{\operatorname{minimize}}} \sum_{t\in \mathcal{T}}^{} p_{{p},t}^{\text{g}+} - \lambda p_{{p},t}^{\text{g}-},
\end{equation}
where $\lambda$ is the applicable feed-in price ratio. In our research, we assumed $\lambda = 0$, which corresponds to maximization of self-consumption.

The prosumer sub-problem is subject to the following constraints:

\textbf{Prosumer power balance:} 
Electrical consumption of prosumer $p$, consisting of grid feed-in power, $p_{{p},t}^{\text{g}-}$, underlying consumption, ${p}_{{p},t}^{\text{}}$, and battery charging power, $ p_{{p},t}^{\text{b}}$, is equal to the power taken from the grid, $p_{{p},t}^{\text{g}+}$, plus the power generated by the PV system, $p_{{p},t}^{\text{pv}}$:
\begin{equation}
p_{{p},t}^{\text{g}+} + p_{{p},t}^{\text{pv}} =  p_{{p},t}^{\text{g}-} +  {p}_{{p},t}^{\text{}} + p_{{p},t}^{\text{b}}.
\end{equation}

\textbf{Battery charge/discharge limits:} 
Battery power should not exceed the charge/discharge limits:
\begin{equation}
\overline{p}_{{p}}^{\text{b}-} \leq {p}_{{p},t}^{\text{b}}  \leq \overline{p}_{{p}}^{\text{b}+},
\end{equation}	
where $\overline{{p}}_{b}^-$ and $\overline{{p}}_{b}^+$ represent the maximum power discharge and charge rates of the prosumer's battery, respectively.		

\textbf{Battery storage capacity limits:} 
Energy stored in a battery of prosumer $p$ should always be less than its capacity:
\begin{equation}
\underline{e}_{{p}}^{\text{b}} \leq {e}_{{p},t}^{\text{b}} \leq \overline{e}_{{p}}^{\text{b}}.
\end{equation}

\textbf{Battery SOC limits:}
Battery SOC is the sum of the  power inflow and the SOC in the previous period:
\begin{align}
& {e}_{{p},t}^{\text{b}} = \eta_p^\text{b} {e}_{{p},t}^{\text{b}} +  {p}_{{p},t}^{\text{b}},  & t \neq 1, \\
& {e}_{{p},t}^{\text{b}} = \eta_p^\text{b} \hat{e}_{{p}}^{\text{b}}  +  {p}_{{p},t}^{\text{b}}, & t=1, 
\end{align}		
where $\hat{e}_{{p}}^{\text{b}}$ represents the initial SOC and is used to establish the connection between adjacent sub-horizons.

%%%\vspace{-1.0em}
\section{Simulation Setup}
The case studies provided in this section compare the computational efficiency of the proposed MST with alternative formulations. For detailed studies on the impact of different technologies on future grids, an interested reader can refer to our previous work~\cite{Riaz15,Riaz16b,Marzooghi16b,Ahmad16a,Marzooghi16,Riaz16a}.

%\vspace{-1.0em}
\subsection{Test System}
We use a modified 14-generator IEEE test system that was initially proposed in \cite{Gibbard10} as a test bed for small-signal analysis. The system is loosely based on the Australian National Electricity Market (NEM), the interconnection on the Australian eastern seaboard. The network is stringy, with large transmission distances and loads concentrated in a few load centres. Generation, demand and the transmission network were modified to meet future load requirements. The modified model consists of 79 buses grouped into four regions, 101 units installed at 14 generation plants and 810 transmission lines.

%\vspace{-1.0em}
\subsection{Test Cases}
To expose the limitations of the different UC formulations, we have selected a typical week with sufficiently varying operating conditions.
Four diverse test cases with different RES penetrations are considered. 
First, RES0 considers only conventional generation, including hydro, black coal, brown coal, combined cycle gas and open cycle gas. The generation mix consists of \SI{2.31}{\giga\watt} hydro, \SI{39.35}{\giga\watt} of coal and \SI{5.16}{\giga\watt} of gas, with the peak load of \SI{36.5}{\giga\watt}. To cater for demand and generation variations, \SI{10}{\percent} reserves are maintained at all times. The generators are assumed to bid at their respective short run marginal costs, based on regional fuel prices~\cite{Tasman09}.

Cases RES30, RES50, RES75 consider, respectively, \SI{30}{\percent}, \SI{50}{\percent} and \SI{75}{\percent} annual energy RES penetration, supplied by wind, PV and CST. Normalized power traces for PV, CST and wind farms (WFs) for the 16-zones of the NEM are taken from the AEMO's planning document~\cite{AEMO2016b}. The locations of RESs are loosely based on the AEMO's 100\% RES study \cite{AEMORES}. 

%\vspace{-1.0em}
\subsection{Modeling Assumptions}
Power traces of all PV modules and wind turbines at one plant are aggregated and represented by a single generator. This is a reasonable assumption given that PV and WF don't provide active power reserves, and are not limited by ramp rates, MUDT, and startup and shutdown costs, which renders the information on the number of online units unnecessary.

Also worth mentioning is that RES can be modeled as negative demand, which can lead to an infeasible solution. Modeling RES (wind and solar PV) as negative demand is namely identical to preventing RES from spilling energy. Given the high RES penetration in future grids, we model RES explicitly as individual generators. 
Unlike solar PV and wind, CST requires a different modeling approach. Given that CST is synchronous generation it also contributes to spinning reserves and system inertia. Therefore, the number of online units in a CST plant needs to be modeled explicitly.

An optimality gap of \SI{1}{}\% was used for all test cases. Simulation were run on Dell OPTIPLEX 9020 desktop computer with Intel(R) Core(TM) i7-4770 CPU with \SI{3.40}{\giga\hertz} clock speed and \SI{16}{\giga B} RAM. 

%%\vspace{-1.0em}
\section{Results and Discussion}
To showcase the computational efficiency of the proposed MST, we first benchmark its performance for different horizon lengths against the BUC formulation employing three binary variables \emph{per unit per time slot} and the AGG formulation where identical units at each plant are aggregated into a single unit, which requires three binary variables \emph{per plant per time slot}. 
We pay particular attention to the techniques used for computational speedup, namely unit clustering, rolling horizon, and constraint clipping. Last, we compare the results of the proposed MST with BUC and AGG formulations for voltage and frequency stability studies.

%\vspace{-1.0em}
\subsection{Binary Unit Commitment (BUC)} \label{BUC}
\begin{figure}
	\centering
	\includegraphics[width=85mm] {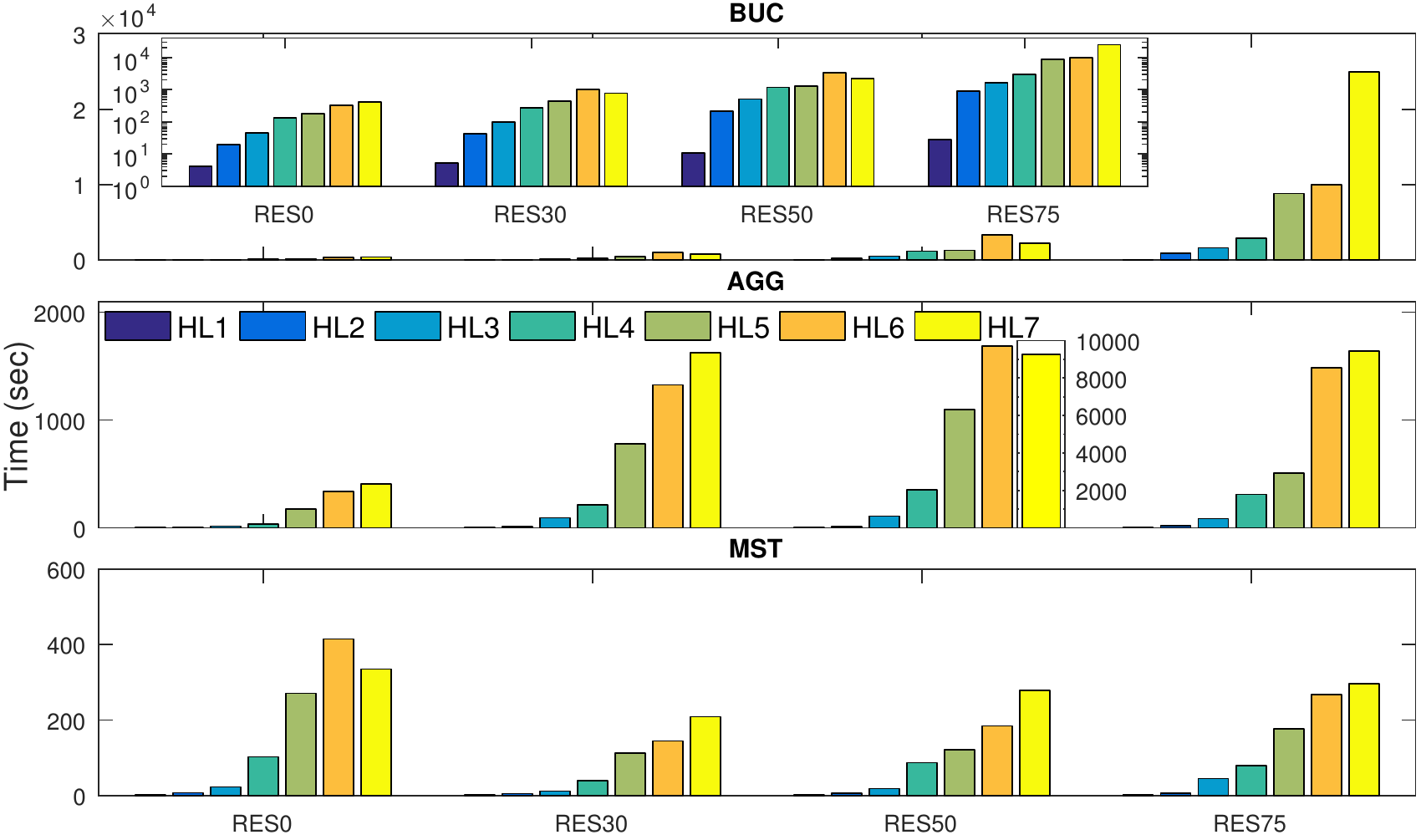}
	%\vspace{-1.0em}
	\caption{Computation time requirements of BUC, AGG and MST for horizon lengths of one (HL1) to seven days (HL7) and different RES penetration levels. In the insets, a logarithmic scale is used for BUC, and a linear scale for AGG case RES50-HL7.}
	\label{fig:BUC_IUC_MST_RT}
%	%\vspace{-1.5em}
\end{figure}

We first run the BUC for horizon lengths varying from one to seven days, Fig.~\ref{fig:BUC_IUC_MST_RT} (top).
As expected, with the increase in the horizon length, the solution time increases exponentially. For a seven-day horizon, the solution time is as high as \SI{25000}{\second} (\SI{7}{\hour}). Observe how the computational burden is highly dependent on the RES penetration. The variability of the RES results in an increased cycling of the conventional thermal fleet, which increases the number of on/off decisions and, consequently the computational burden. In addition to that, a higher RES penetration involves an increased operation of CST. This poses an additional computational burden due to the decision variables associated with TES that span several time slots.
In summary, the computational burden of the BUC renders it inappropriate for scenario analysis involving extended horizons.

%\vspace{-1.0em}
\subsection{Aggregated Formulation (AGG)}
Aggregating identical units at a power plant into a single unit results in a smaller number of binary variables, which should in principle reduce the computational complexity. 
Fig.~\ref{fig:BUC_IUC_MST_RT} confirms that this is mostly true, however, for RES50-HL7 the computation time is higher than in the BUC formulation. The reason for that is that, in this particular case, the BUC formulation has a tighter relaxation than the AGG formulation and, consequently, a smaller root node gap. Compared to the MST formulation, with a similar number of variables than the AGG formulation, the MST has considerably shorter computation time due to a smaller root node gap. 
 
\begin{figure}
	\centering
	\includegraphics[width=85mm] {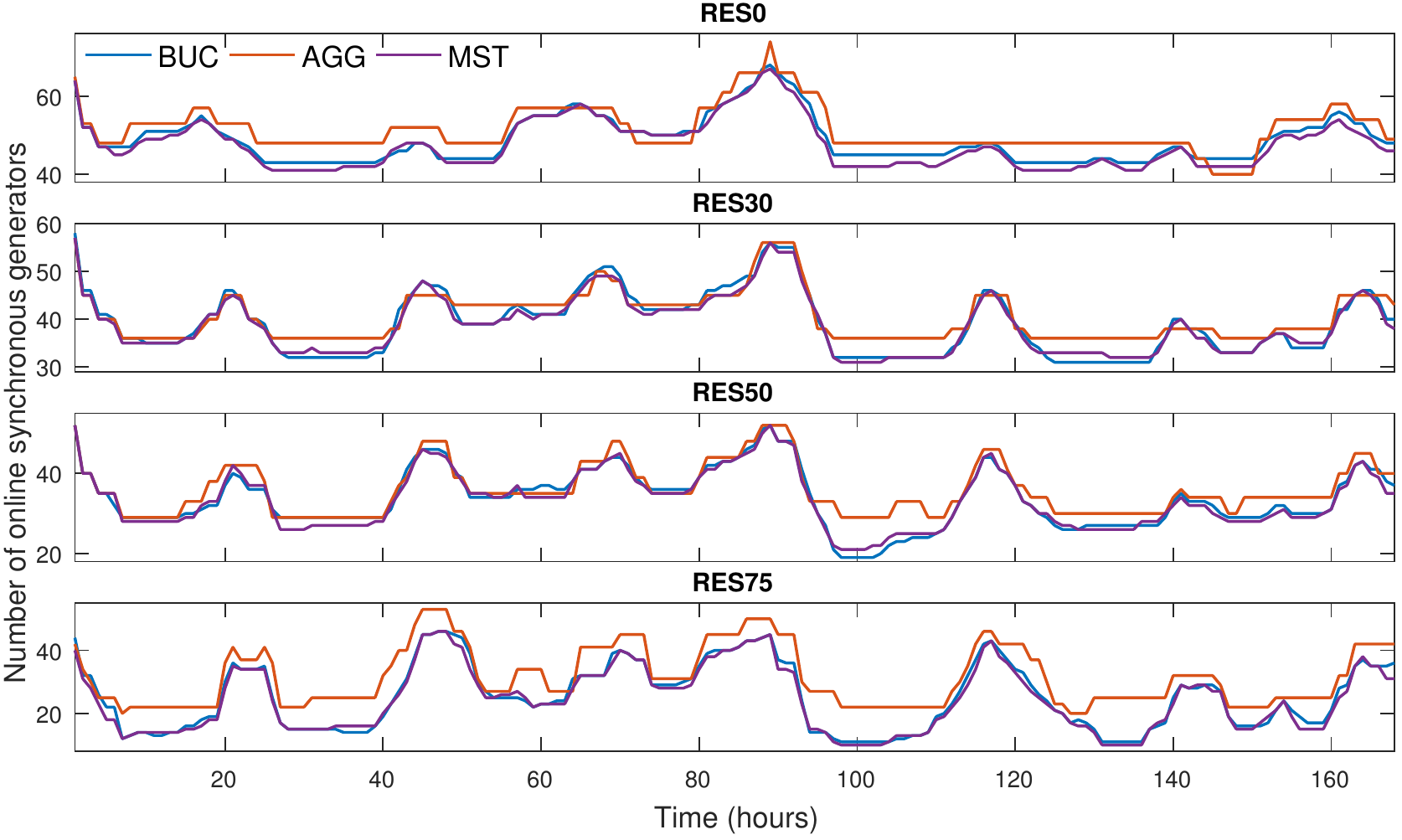}
	%\vspace{-1.0em}
	\caption{Total number of online synchronous generators for BUC, AGG and MST for different RES penetration levels and a horizon length of seven days.}
	\label{fig:online_gen}
%	%\vspace{-1.5em}
\end{figure}

In terms of accuracy, the AGG formulation works well for balancing studies~\cite{Riaz15,Riaz16a}. On the other hand, the number of online synchronous generators in the dispatch differs significantly from the BUC, which negatively affects the accuracy of voltage and frequency stability analysis, as shown later. Due to a large number of online units in a particular scenario, a direct comparison of dispatch levels and reserves from each generator is difficult. Therefore, we compare the total number of online synchronous generators, which serves as a proxy to the available system inertia. Fig.~\ref{fig:online_gen} shows the number of online generators of four different RES penetration levels for a horizon length of seven days. For most of the hours there is a significant difference between the number of online units obtained from the BUC and the AGG formulation.     

In conclusion, despite its computational advantages, the AGG formulation is not appropriate for stability studies due to large variations in the number of online synchronous units in the dispatch results. In addition to that, the computational time is comparable to the BUC in some cases.

%\vspace{-1.0em}\subsection{Computational Speedup Assessment}
We now evaluate the effectiveness of the techniques for the computational speedup.

\subsubsection{Unit Clustering}
In unit clustering, binary variables associated with the generation unit constraints are replaced with a smaller number of integer variables, which allows aggregating several identical units into one equivalent unit, but with the number of online units retained. This results in a significant reduction in the number of variables and, consequently, in the computational speedup. Compared to the BUC, the number of variables in the MST with this technique alone reduces from \SI{24649}{} to \SI{5990}{} for RES75 with a horizon length of seven days. Therefore, the solution time for RES75-HL7 reduces from \SI{25000}{\second} in the BUC to \SI{450}{\second} in MST with unit clustering alone.

\subsubsection{Rolling Horizon Approach}
A rolling horizon approach splits the UC problem into shorter horizons. Given the exponential relationship between the computational burden and the horizon length, as discussed in Section~\ref{BUC}, solving the problem in a number of smaller chunks instead of in one block results in a significant computational speedup. The accuracy and the consistency of the solution are maintained by having an appropriate overlap between the adjacent horizons. However, the overlap depends on the time constants of the problem. Long term storage, for example, might require longer solution horizons. The solution times for different RES penetrations are shown in Table~\ref{table:Time_comp}. Observe that in the RES75 case, the effect of rolling horizon is much more pronounced, which confirms the validity of the approach for studies with high RES penetration.
%\begin{figure}
%	%%%\vspace{-1.0em}
%	\centering
%	\includegraphics[width=85mm] {Fig_On-Line_RH_V2.pdf}
%	%%\vspace{-0.3em}
%	\caption{Comparison of the number of online synchronous generators using BUC for RES50 with a rolling horizon (RH) approach for sub-horizon lengths of two to six days with one day overlap against a seven day long monolithic block (MB).}
%	\label{fig:RH_eff_BUC}
%	%%%\vspace{-1.5em}
%\end{figure}
\begin{table}[!t]
\renewcommand{\arraystretch}{1.0}
\centering
\caption{Computation time comparison (MB = monolithic block, RH = rolling horizon, 7 = 7 days, 2+1 = 2 days with one day overlap).}
\label{table:Time_comp}
	\begin{tabular}{|c|c|c|c|c|}
		\hline            & RES0      & RES30     & RES50     & RES75 \\                           
		                  & (minutes) & (minutes) & (minutes) & (minutes) \\
		\hline BUC MB 7   & 6.92      & 12.95     & 37.11     & 415.25 \\
		\hline AGG MB 7   &   6.81    &    27.08   &   154.27    & 27.37 \\
		\hline MST MB 7   & 2.12      & 3.34      & 4.73      & 5.32  \\
		\hline BUC RH 2+1 & 2.38      & 4.03      & 24.18     & 74.70\\
		\hline AGG RH 2+1 &  0.15     &   0.20    &   0.27    & 0.25 \\
		\hline MST RH 2+1 & 0.35      & 0.71      & 0.60      & 0.76\\
		\hline
	\end{tabular}
%	%\vspace{-1.0em}
\end{table}

\subsubsection{Constraint Clipping}
Eliminating non binding constraints can speedup the computation even further. Table~\ref{table:CR_eff} shows the number of constraints for different scenarios with and without constraint clipping. Observe that the number of redundant constraints is higher in scenarios with a higher RES penetration. The reason is that a higher RES penetration requires more flexible gas generation with ramp rates shorter than the time resolution (one hour in our case). Note that the benefit of constraint clipping with a shorter time resolution will be smaller.
\begin{table}
	\renewcommand{\arraystretch}{1.0}
	\centering
	\caption{The impact of constraint clipping (CC) on the total number of constraints for all cases with a horizon length of seven days.}
	\label{table:CR_eff}
	\begin{tabular}{|c|c|c|c|}%{|c|c|}
	\hline Case & \multicolumn{3}{c|}{ Number of Constraints}\\		
	 \cline{2-4} & without CC&with CC & \% reduction\\
	 \hline RES0  & 64555 & 61332 & 4.99\\
	 \hline RES30 & 62520 & 56617 & 9.44\\
	 \hline RES50 & 62500 & 56777 & 9.15\\
	 \hline RES75 & 62740 & 57017 & 9.12\\
	 \hline
	\end{tabular}
%	%\vspace{-2.0em}
\end{table}

%\vspace{-1.0em}
\subsection{MST Computation Time and Accuracy}  
The proposed MST outperforms the BUC and AGG in terms of the computational time by several orders of magnitude, as shown in Fig.~\ref{fig:BUC_IUC_MST_RT} (bottom). The difference is more pronounced at higher RES penetration levels. For RES75, the MST is more than \SI{500}{} times faster than the BUC. In terms of the accuracy, the MST results are almost indistinguishable from the BUC results, as evident from Fig.~\ref{fig:online_gen} that shows the number of online synchronous units for different RES penetration levels. Minor differences in the results stem from the nature of the optimization problem. Due to its mixed-integer structure, the problem is non-convex and has therefore several local optima. Given that the BUC and the MST are mathematically not equivalent, the respective solutions might not be exactly the same. The results are nevertheless very close, which confirms the validity of the approach for the purpose of scenario analysis. The loadability and inertia results presented later further support this conclusion.
 
%\vspace{-1.0em}
\subsection{Stability Assessment}
To showcase the applicability of the MST for stability assessment, we analyze system inertia and loadability that serve as a proxy to frequency and voltage stability, respectively. More detailed stability studies are covered in our previous work, including small-signal stability \cite{Liu2016}, frequency stability \cite{Ahmad16a}, and voltage stability \cite{Riaz16b}.

\subsubsection{System inertia}
Fig.~\ref{fig:Voltage_Stability} (bottom) shows the system inertia for the BUC, AGG and the proposed MST, respectively, for RES0. Given that the inertia is the dominant factor in the frequency response of a system after a major disturbance, the minuscule difference between the BUC and the MST observed in Fig.~\ref{fig:Voltage_Stability} validates the suitability of the MST for frequency stability assessment. The inertia captured by the AGG, on the other hand, is either over or under estimated and so does not provide a reliable basis for frequency stability assessment.

\subsubsection{Loadability Analysis}
\begin{figure}
	\centering
	\includegraphics[width=85mm] {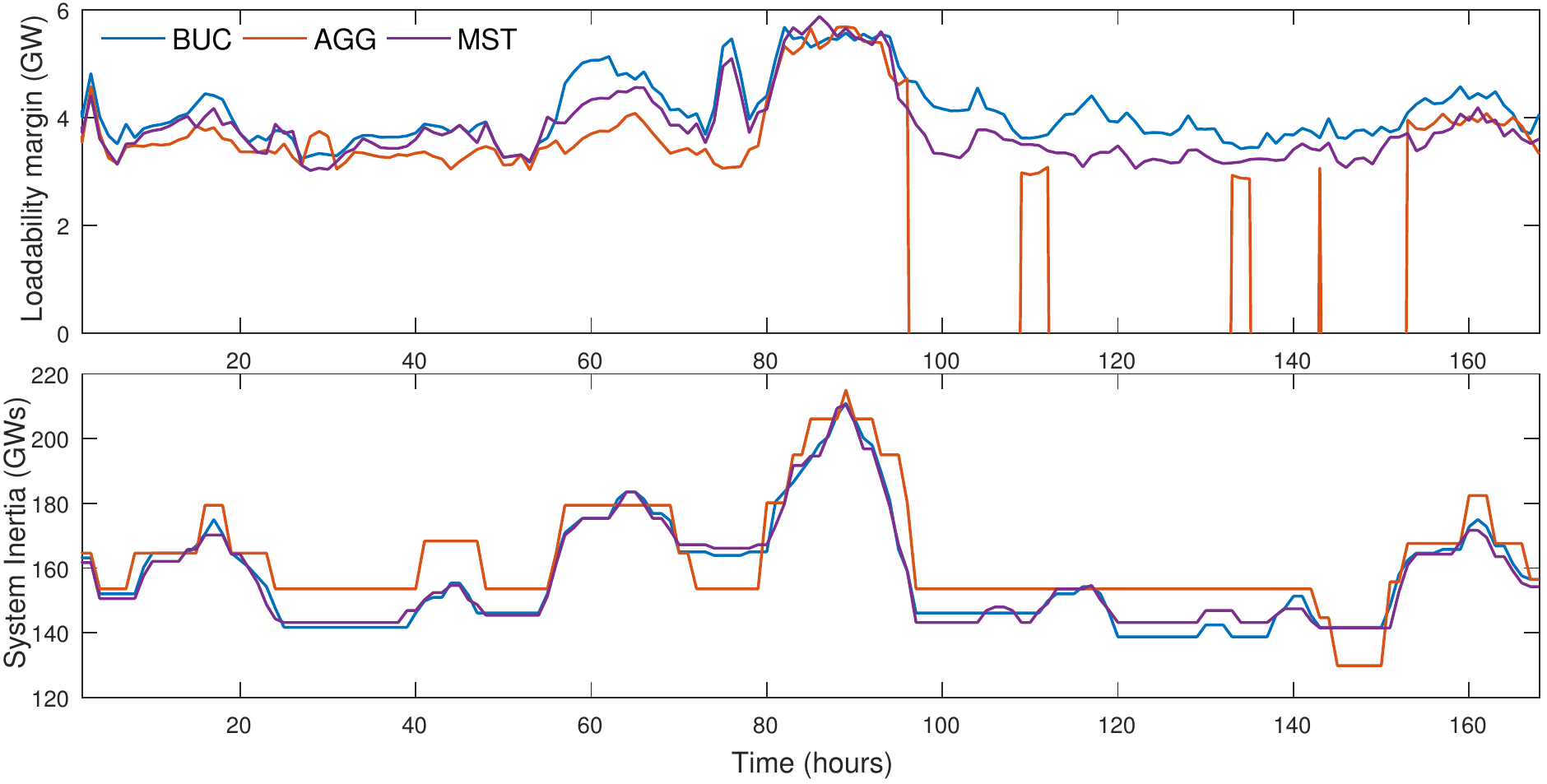}
	%\vspace{-1em}
	\caption{Loadability margins (top) and system inertia (bottom) computed based on dispatch results of different techniques i.e BUC, AGG and proposed MST for RES0.}
	\label{fig:Voltage_Stability}
%	%\vspace{-1.25em}
\end{figure}

The dispatch results from the MST are used to calculate power flows, which are then used in loadability analysis\footnote{The loadability analysis is performed by uniformly increasing the load in the system until the load flow fails to converge. The loadability margin is calculated as the difference between the base system load and the load in the last convergent load flow iteration.}. Fig.~\ref{fig:Voltage_Stability} (top) shows loadability margins for the RES0 scenario for different UC formulations. Observe that the BUC and the MST produce very similar results. The AGG formulation, on the other hand, gives significantly different results. From hours \SI{95}{} to \SI{150}{}, in particular, the AGG results show that the system is unstable most of the time, which is in direct contradiction to the accurate BUC formulation. 
Compared to the inertia analysis, the differences between the formulations are much more pronounced.
Unlike voltage, frequency is a system variable, which means that it is uniform across the system. In addition to that, inertia only depends on the number of online units but not on their dispatch levels.
Voltage stability, on the other hand, is highly sensitive both to the number of online units and their dispatch levels, which affects the available reactive power support capability, as illustrated in Fig. \ref{fig:Generation_comparision}. 
Close to the voltage stability limit, the system becomes highly nonlinear, so even small variations in dispatch results can significantly change the power flows and, consequently, voltage stability of the system. One can argue that in comparison to BUC  the proposed MST result in the more conservative loadability margin, although this is not always the case (around hour \SI{85}{}, the MST is less conservative).
%\vspace{-1em}

%%\vspace{-1.0em}
\section{Conclusion}
This paper has proposed a computationally efficient electricity market simulation tool based on a UC problem suitable for future grid scenario analysis. The proposed UC formulation includes an explicit network representation and accounts for the uptake of emerging demand side technologies in a unified generic framework while allowing for a subsequent stability assessment. We have shown that unit aggregation, used in conventional planning-type UC formulations to achieve computational speedup, fails to properly capture the system inertia and reactive power support capability, which is crucial for stability assessment. To address this shortcoming, we have proposed a UC formulation that models the number of online generation units explicitly and is amenable to a computationally expensive time-series analysis required in future grid scenario analysis. To achieve further speedup, we use a rolling horizon approach and constraint clipping.

The effectiveness of the computational speedup techniques depends on the problem structure and the technologies involved so the results cannot be readily generalized. The computational speedup varies between \SI{20}{} to more than \SI{500}{} times, for a zero and \SI{75}{}\% RES penetration, respectively, which can be explained by a more frequent cycling of the conventional thermal units in the high-RES case. The simulation results have shown that the computational speedup doesn't jeopardize the accuracy. Both the number of online units that serves as a proxy for the system inertia and the loadability results are in close agreement with more detailed UC formulations, which confirms the validity of the approach for long term future grid studies, where one is more interested in finding weak points in the system rather than in a detailed analysis of an individual operating condition.

\ifCLASSOPTIONcaptionsoff
  \newpage
\fi

%%\vspace{-1.75em}

\bibliographystyle{IEEEtran}
\bibliography{UC}
\end{document}